\documentclass{article}
\usepackage{amssymb}

\input{tcilatex}
\begin{document}

\begin{center}
\bigskip \textbf{Multipoint Cauchy problem for Schr\"{o}dinger type
equations with general elliptic part }

\textbf{Veli Shakhmurov}

Department of Mechanical Engineering, Okan University, Akfirat, Tuzla 34959
Istanbul, Turkey,

E-mail: veli.sahmurov@okan.edu.tr

\textbf{Abstract}
\end{center}

In this paper, the existence, uniqueness and regularity properties,
Strichartz type estimates for solution of multipoint Cauchy problem for
linear and nonlinear Schr\"{o}dinger equations with general elliptic leading
part is obtained.

\textbf{Key Word:}$\mathbb{\ \ }$Schr\H{o}dinger equations\textbf{, }%
elliptic operators\textbf{, }Semigroups of operators, local solutions

\begin{center}
\ \ \textbf{AMS 2010: 35Q41, 35K15, 47B25, 47Dxx, 46E40 }

\textbf{1. Introduction}
\end{center}

\bigskip\ Consider the multipoint Cauchy problem for nonlinear Schr\"{o}%
dinger equations (NLS)

\begin{equation}
i\partial _{t}u+Lu+F\left( u\right) =0,\text{ }x\in R^{n},\text{ }t\in \left[
0,T\right] ,  \tag{1.1}
\end{equation}%
\begin{equation}
u\left( 0,x\right) =\varphi \left( x\right) +\dsum\limits_{k=1}^{m}\alpha
_{k}u\left( \lambda _{k},x\right) ,\text{ for a.e. }x\in R^{n},  \tag{1.2}
\end{equation}%
where $L$ is an elliptic operator defined by

\begin{equation}
Lu=\dsum\limits_{i,j=1}^{2}a_{ij}\frac{\partial ^{2}u}{\partial
x_{i}\partial x_{j}},\text{ }a_{ij}\in \mathbb{C},  \tag{1.3}
\end{equation}%
$m$ is a positive integer, $\alpha _{k}$ are complex numbers, $\lambda
_{k}\in \left( 0,\right. \left. T\right] ,$ $F$ is a nonlinear operator, $%
\mathbb{C}-$denotes the set of complex numbers and $u=$ $u(t,x)$ is the
unknown function. If $F\left( u\right) =$ $\lambda \left\vert u\right\vert
^{p}u$ in $\left( 1.1\right) $ we get the multipoint Cauchy problem
nonlinear equation%
\begin{equation}
i\partial _{t}u+Lu+\lambda \left\vert u\right\vert ^{p}u=0,\text{ }x\in
R^{n},\text{ }t\in \left[ 0,T\right] ,  \tag{1.4}
\end{equation}%
\[
u\left( 0,x\right) =\varphi \left( x\right) +\dsum\limits_{k=1}^{m}\alpha
_{k}u\left( \lambda _{k},x\right) \text{ for a.e. }x\in R^{n}, 
\]%
where $p\in \left( 1,\infty \right) $, $\lambda $ is a real number.

\bigskip By rescaling the values of $u$\ it is possible to restrict
attention to the cases $\lambda =1$ or $\lambda =-1.$ These call as the
focusing and defocusing Schr\"{o}dinger equations, respectively. The
equation $\left( 1.1\right) $ also contain two critical case. These are the
mass-critical Schr\"{o}dinger equation,%
\[
i\partial _{t}u+Lu+\lambda \left\vert u\right\vert ^{\frac{4}{n}}=0,\text{ }%
x\in R^{n},\text{ }t\in \left[ 0,T\right] ,
\]%
which is associated with the conservation of mass,

\[
M\left( u\left( t\right) \right) :=\dint\limits_{R^{n}}\left\Vert u\left(
t,x\right) \right\Vert _{E}^{2}dx 
\]%
and the energy-critical Schr\"{o}dinger equation (in dimensions $n>2$),

\begin{equation}
i\partial _{t}u+Lu+\lambda \left\vert u\right\vert ^{\frac{4}{n-2}}=0,\text{ 
}x\in R^{n},\text{ }t\in \left[ 0,T\right] ,  \tag{1.5}
\end{equation}%
which is associated with the conservation of energy,%
\[
H\left( u\left( t\right) \right) :=\dint\limits_{R^{n}}\left[ \frac{1}{2}%
\left\vert \left( Lu,u\right) \left( t,x\right) \right\vert ^{2}+\lambda
\left( \frac{1}{2}-\frac{1}{n}\right) \left\vert u\left( t,x\right)
\right\vert ^{\frac{2n}{n-2}}\right] dx, 
\]%
where $\left( Lu,u\right) $ denotes scalar product of $Lu$ and $u$ in $%
L^{2}\left( R^{n}\right) .$

The existence of solutions and regularity properties of Cauchy problem for
NLS equations studied e.g in $\left[ 2-10\right] $, $\left[ \text{ }14,\text{
}16\right] $ and the references therein.\ In contrast, to the mentioned
above results we will study the existence, uniqueness and the regularity
properties of the multipoint Cauchy problem $\left( 1.1\right) -\left(
1.2\right) $.

\begin{center}
\textbf{2. Definitions and background}
\end{center}

Let $L_{t}^{q}L_{x}^{r}\left( \left( a,b\right) \times \Omega \right) $
denotes the space of strongly measurable functions that are defined on the
measurable set $\left( a,b\right) \times \Omega $ with the norm 
\[
\left\Vert f\right\Vert _{L_{t}^{q}L_{x}^{r}\left( \left( a,b\right) \times
\Omega \right) }=\left( \dint\limits_{a}^{b}\left( \int\limits_{\Omega
}\left\vert f\left( t,x\right) \right\vert ^{r}dx\right) ^{\frac{q}{r}%
}dt\right) ^{\frac{1}{q}},\text{ }1\leq q,r<\infty \ . 
\]

Let $\digamma $ denotes the Fourier transformation, $\hat{u}=\digamma u$ and 
\[
s\in \mathbb{R},\text{ }\xi =\left( \xi _{1},\xi _{2},...,\xi _{n}\right)
\in R^{n},\text{ }\left\vert \xi \right\vert ^{2}=\dsum\limits_{k=1}^{n}\xi
_{k}^{2}, 
\]%
\[
\langle \xi \rangle =\left( 1+\left\vert \xi \right\vert ^{2}\right) ^{\frac{%
1}{2}}. 
\]%
$S=S(R^{n})$ denotes the Schwartz class, i.e. the space of all
complex-valued rapidly decreasing smooth functions on $R^{n}$ equipped with
its usual topology generated by seminorms. $S(R^{n})$ denoted by just $S$.
Let $S^{\prime }(R^{n})$ denote the space of all continuous linear
operators, $L:S\rightarrow \mathbb{C}$, equipped with the bounded
convergence topology. Recall $S(R^{n})$ is norm dense in $L^{p}(R^{n})$ when 
$1<p<\infty .$ Let $D^{\prime }\left( \Omega \right) $ denote the class of
generalized functions on $\Omega \subset R^{n}$. Consider Sobolev space $%
W^{s,p}\left( R^{n}\right) $ and homogeneous Sobolev spaces $\mathring{W}%
^{s,p}\left( R^{n}\right) $\ defined by respectively, 
\[
W^{s,p}\left( R^{n}\right) =\left\{ u:u\in S^{\prime }(R^{n}),\right. \text{ 
} 
\]%
\[
\left\Vert u\right\Vert _{W^{s,p}\left( R^{n}\right) }=\left. \left\Vert
\digamma ^{-1}\left( 1+\left\vert \xi \right\vert ^{2}\right) ^{\frac{s}{2}}%
\hat{u}\right\Vert _{L^{p}\left( R^{n}\right) }<\infty \right\} , 
\]

\[
\mathring{W}^{s,p}\left( R^{n}\right) =\left\{ u:u\in S^{\prime
}(R^{n}),\left\Vert u\right\Vert _{\mathring{W}^{s,p}\left( R^{n}\right)
}=\left\Vert \digamma ^{-1}\left\vert \xi \right\vert ^{s}\right\Vert
_{L^{p}\left( R^{n}\right) }<\infty \right\} .\text{ } 
\]

Sometimes we use one and the same symbol $C$ without distinction in order to
denote positive constants which may differ from each other even in a single
context. When we want to specify the dependence of such a constant on a
parameter, say $\alpha $, we write $C_{\alpha }$.

Let $\ L$ is differential operator defined by $\left( 1.4\right) .$

\textbf{Condition 2.1. }Assume $a_{ij}=a_{ji}$ and there are positive
constants $M_{1}$ and $M_{2}$ such that $\ M_{1}\left\vert \xi \right\vert
^{2}\leq L\left( \xi \right) \leq M_{2}\left\vert \xi \right\vert ^{2}$ for $%
\xi =\left( \xi _{1},\xi _{2},...\xi _{n}\right) \in R^{n}$, where%
\[
\left\vert \xi \right\vert ^{2}=\dsum\limits_{k=1}^{n}\xi _{k}^{2}\text{, }%
L\left( \xi \right) =\dsum\limits_{i,j=1}^{2}a_{ij}\xi _{i}\xi _{j}. 
\]

\textbf{Definition 2.2.} Consider the initial value problem $(1.1)-\left(
1.2\right) $ for $\varphi \in \mathring{W}^{s,p}\left( R^{n}\right) $. This
problem is critical when $s=s_{c}:=\frac{n}{2}-\frac{2}{p},$ subcritical
when $s>s_{c}$, and supercritical when $s<s_{c}$.

We write $a\lesssim b$ to indicate that $a\leq Cb$ for some constant $C$,
which is permitted to depend on some parameters.\bigskip

\begin{center}
\textbf{3. Dispersive and Strichartz type inequalities} \textbf{for linear
Schr\"{o}dinger equation}
\end{center}

\ Let the operator $iL$ generates a continious $C_{0}$ group $e^{itL\left(
\xi \right) }.$ It can be shown that the fundamental solution of the free
Schr\"{o}dinger equation%
\begin{equation}
i\partial _{t}u+Lu=0,\text{ }t\in \left[ 0,T\right] ,\text{ }x\in R^{n} 
\tag{3.1}
\end{equation}%
can be expressed as $U_{L}\left( t\right) =\digamma ^{-1}e^{itL\left( \xi
\right) }$, i.e.

\begin{equation}
U_{L}\left( t\right) f\left( x\right) =\dint\limits_{R^{n}}U_{L}\left(
t\right) \left( x-y\right) f\left( y\right) dy.\text{ }  \tag{3.2}
\end{equation}

\textbf{Lemma 3.1.}The following dispersive inequalities hold 
\begin{equation}
\left\Vert U_{L}\left( t\right) f\right\Vert _{L_{x}^{p}\left( R^{n}\right)
}\lesssim t^{-\left[ n\left( \frac{1}{2}-\frac{1}{p}\right) \right]
}\left\Vert f\right\Vert _{L_{x}^{p^{\prime }}\left( R^{n}\right) }, 
\tag{3.3}
\end{equation}%
\begin{equation}
\left\Vert U_{L}\left( t-s\right) f\right\Vert _{L^{\infty }\left(
R^{n}\right) }\lesssim \left\vert t-s\right\vert ^{-\frac{n}{2}}\left\Vert
f\right\Vert _{L^{1}\left( R^{n}\right) }  \tag{3.4}
\end{equation}%
for $t\neq 0,$ $2\leq p\leq \infty ,$ $\frac{1}{p}+\frac{1}{p^{\prime }}=1.$

\textbf{Proof. }Indeed, by using Young's integral inequality from $\left(
3.2\right) $ we get 
\begin{equation}
\left\Vert U_{L}\left( t\right) f\right\Vert _{L_{x}^{p}\left( R^{n}\right)
}\lesssim \left\vert t\right\vert ^{-n\left( \frac{1}{2}-\frac{1}{p}\right)
}\left\Vert f\right\Vert _{L_{x}^{p^{\prime }}\left( R^{n}\right) }, 
\tag{3.5}
\end{equation}%
\begin{equation}
\left\Vert U_{L}\left( t\right) f\right\Vert _{L_{x}^{\infty }\left(
R^{n}\right) }\lesssim \left\vert t\right\vert ^{-\frac{n}{2}}\left\Vert
f\right\Vert _{L_{x}^{1}\left( R^{n}\right) }.  \tag{3.6}
\end{equation}

\textbf{Condition 3.1. }Assume $n\geq 1,$%
\[
\frac{2}{q}+\frac{n}{r}\leq \frac{n}{2},\text{ }2\leq q,r\leq \infty \text{
\ and }\left( n,\text{ }q,\text{ }r\right) \neq \left( 2,\text{ }2,\text{ }%
\infty \right) . 
\]

\textbf{Remark 3.1. }If $\frac{2}{q}+\frac{n}{r}=\frac{n}{2},$ then $(q,$ $%
r) $ is called sharp admissible, otherwise $(q,$ $r)$ is called nonsharp
admissible. Note in particular that when $n>2$ the endpoint $\left( 2\text{, 
}\frac{2n}{n--2}\right) $ is called sharp admissible.

\bigskip For a space-time slab $\left[ 0,T\right] \times R^{n}$, we define
the $E-$valued Strichartz norm%
\[
\left\Vert u\right\Vert _{S^{0}\left( I\right) }=\sup\limits_{\left(
q,r\right) \text{ admissible}}\left\Vert u\right\Vert
_{L_{t}^{q}L_{x}^{r}\left( I\times R^{n}\right) }, 
\]%
where $S^{0}\left( \left[ 0,T\right] \right) $ is the closure of test
functions under this norm and $N^{0}\left( \left[ 0,T\right] \right) $
denotes the dual of $S^{0}\left( \left[ 0,T\right] \right) .$

Assume $H$ is an abstract Hilbert space and $Q$ is a Hilbert space of
function. Suppose for each $t\in \mathbb{R}$ an operator $U\left( t\right) $%
: $Q\rightarrow L^{2}\left( \Omega \right) $ obeys the following estimates:

\begin{equation}
\left\Vert U\left( t\right) f\right\Vert _{L_{x}^{2}\left( \Omega \right)
}\lesssim \left\Vert f\right\Vert _{H}  \tag{3.7}
\end{equation}%
for all $t,$ $\Omega \subset R^{n}$ and all $f\in Q;$%
\begin{equation}
\left\Vert U\left( s\right) U^{\ast }\left( t\right) g\right\Vert
_{L_{x}^{\infty }\left( \Omega \right) }\lesssim \left\vert t-s\right\vert
^{-\frac{n}{2}}\left\Vert g\right\Vert _{L_{x}^{1}\left( \Omega \right) } 
\tag{3.8}
\end{equation}%
\begin{equation}
\left\Vert U\left( s\right) U^{\ast }\left( t\right) g\right\Vert
_{L_{x}^{\infty }\left( \Omega \right) }\lesssim \left( 1+\left\vert
t-s\right\vert ^{-\frac{n}{2}}\right) \left\Vert g\right\Vert
_{L_{x}^{1}\left( \Omega \right) }  \tag{3.9}
\end{equation}%
for all $t\neq s$ and all $g\in L_{x}^{1}\left( \Omega \right) .$

For proving the main theorem of this section, we will use the following
bilinear interpolation result (see $\left[ 1\right] $, Section 3.13.5(b)).

\bigskip \textbf{Lemma 3.2. }Assume $A_{0}$, $A_{1},$ $B_{0}$, $B_{1},$ $%
C_{0}$, $C_{1}$ are Banach spaces and $T$ is a bilinear operator bounded
from ($A_{0}\times B_{0}$, $A_{0}\times B_{1},$ $A_{1}\times B_{0}$ ) into ($%
C_{0}$, $C_{1}$, $C_{1}$), respectively. Then whenever $0<\theta _{0},$ $%
\theta _{1}<\theta <1$ are such that $1\leq \frac{1}{p}+\frac{1}{q}$ and $%
\theta =\theta _{0}+$ $\theta _{1}$, the operator is bounded from 
\[
\left( A_{0}\text{, }A_{1}\right) _{\theta _{0}pr}\times \left( B_{0}\text{, 
}B_{1}\right) _{\theta _{1}qr} 
\]%
to $\left( C_{0}\text{, }C_{1}\right) _{\theta r}.$

By following $\left[ \text{9, Theorem 1.2}\right] $ we have:

\textbf{Theorem 3.1. }Assume $U(t)$ obeys $\left( 3.8\right) $ and $\left(
3.9\right) $. Then the following estimates are hold%
\begin{equation}
\left\Vert U\left( t\right) f\right\Vert _{L_{t}^{q}L_{x}^{r}}\lesssim
\left\Vert f\right\Vert _{H},  \tag{3.10}
\end{equation}

\begin{equation}
\left\Vert \dint U^{\ast }\left( s\right) F\left( s\right) ds\right\Vert
_{Q}\lesssim \left\Vert F\right\Vert _{L_{t}^{q^{\prime }}L_{x}^{r^{\prime
}}},  \tag{3.11}
\end{equation}%
\begin{equation}
\dint\limits_{s<t}\left\Vert A^{\alpha }U\left( t\right) U^{\ast }\left(
s\right) F\left( s\right) ds\right\Vert _{L_{t}^{q}L_{x}^{r}}\lesssim
\left\Vert F\right\Vert _{L_{t}^{\tilde{q}^{\prime }}L_{x}^{\tilde{r}%
^{\prime }}},  \tag{3.12}
\end{equation}%
for all sharp admissible exponent pairs $\left( q,r\right) $, $\left( \tilde{%
q},\tilde{r}\right) .$ Furthermore, if the decay hypothesis is strengthened
to $(3.9)$, then $(3.10)$, $(3.11)$ and $(3.12)$ hold for all admissible $%
\left( q,\text{ }r\right) $, $\left( \tilde{q},\text{ }\tilde{r}\right) .$

\textbf{Proof. The first step: }Consider the nonendpoint case, i.e. $\left(
q,\text{ }r\right) \neq $ $\left( 2,\text{ }\frac{2n}{n-2}\right) $ and will
show firstly, the estimates $\left( 3.10\right) $, $\left( 3.11\right) .$ By
duality, $(3.10)$ is equivalent to $(3.11)$. By the $TT^{\ast }$ method, $%
(3.11)$ is in turn equivalent to the bilinear form estimate 
\begin{equation}
\left\vert \dint \dint \langle U^{\ast }\left( s\right) F\left( s\right)
,U^{\ast }\left( t\right) G\left( t\right) \rangle dsdt\right\vert \lesssim
\left\Vert F\right\Vert _{L_{t}^{q^{\prime }}L_{x}^{r^{\prime }}}\left\Vert
G\right\Vert _{L_{t}^{q^{\prime }}L_{x}^{r^{\prime }}}.  \tag{3.13}
\end{equation}

By symmetry it suffices to show the to the retarded version of $\left(
3.13\right) $%
\begin{equation}
\left\vert T\left( F,G\right) \right\vert \lesssim \left\Vert F\right\Vert
_{L_{t}^{q^{\prime }}L_{x}^{r^{\prime }}}\left\Vert G\right\Vert
_{L_{t}^{q^{\prime }}L_{x}^{r^{\prime }}},  \tag{3.14}
\end{equation}%
where $T\left( F,G\right) $ is the bilinear form defined by 
\[
T\left( F,G\right) =\dint \dint\limits_{s<t}\langle U\left( s\right) ^{\ast
}F\left( s\right) ,\left( U\left( t\right) \right) ^{\ast }G\left( t\right)
\rangle dsdt 
\]

By real interpolation between the bilinear form of $\left( 3.7\right) $ we
get 
\[
\left\vert \langle \left( U\left( s\right) \right) ^{\ast }F\left( s\right)
,\left( U\left( t\right) \right) ^{\ast }G\left( t\right) \rangle
\right\vert \lesssim \left\Vert F\left( s\right) \right\Vert
_{L_{x}^{2}}\left\Vert G\left( t\right) \right\Vert _{L_{x}^{2}}. 
\]%
By using the bilinear form of $\left( 3.8\right) $ we have 
\begin{equation}
\left\vert \langle \left( U\left( s\right) \right) ^{\ast }F\left( s\right)
,\left( U\left( t\right) \right) ^{\ast }G\left( t\right) \rangle
\right\vert \lesssim  \tag{3.15}
\end{equation}%
\[
\left\vert t-s\right\vert ^{-\frac{n}{2}}\left\Vert F\left( s\right)
\right\Vert _{L_{x}^{1}\left( \Omega \right) }\left\Vert G\left( t\right)
\right\Vert _{L_{x}^{1}\left( \Omega \right) }. 
\]%
In a similar way, we obtain 
\begin{equation}
\left\vert \langle \left( U\left( s\right) \right) ^{\ast }F\left( s\right)
,\left( U\left( t\right) \right) ^{\ast }G\left( t\right) \rangle
\right\vert \lesssim  \tag{3.16}
\end{equation}%
\[
\left\vert t-s\right\vert ^{--1-\beta \left( r,r\right) }\left\Vert F\left(
s\right) \right\Vert _{L_{x}^{r^{\prime }}\left( \Omega \right) }\left\Vert
G\left( t\right) \right\Vert _{L_{x}^{r^{\prime }}\left( \Omega \right) }, 
\]%
where $\beta (r,\tilde{r})$ is given by 
\begin{equation}
\beta (r,\tilde{r})=\frac{n}{2}-1-\frac{n}{2}\left( \frac{1}{r}-\frac{1}{%
\tilde{r}}\right) .  \tag{3.17}
\end{equation}

It is clear that $\beta (r,r)\leq 0$ when $n>2.$ In the sharp admissible
case we have 
\[
\frac{1}{q}+\frac{1}{q^{\prime }}=-\beta (r,r), 
\]%
and $(3.14)$ follows from $(3.16)$ and the Hardy-Littlewood-Sobolev
inequality ($[20]$) when $q>q^{\prime }.$

If we are assuming the truncated decay $(3.9)$, then $(3.16)$ can be
improved to 
\begin{equation}
\left\vert \langle \left( U\left( s\right) \right) ^{\ast }F\left( s\right)
,\left( U\left( t\right) \right) ^{\ast }G\left( t\right) \rangle
\right\vert \lesssim  \tag{3.18}
\end{equation}%
\[
\left( 1+\left\vert t-s\right\vert \right) ^{-1-\beta \left( r,r\right)
}\left\Vert F\left( s\right) \right\Vert _{L_{x}^{r^{\prime }}\left( \Omega
\right) }\left\Vert G\left( t\right) \right\Vert _{L_{x}^{r^{\prime }}\left(
\Omega \right) } 
\]%
and now Young's inequality gives $(3.14)$ when 
\[
-\beta (r,r)+\frac{1}{q}>\frac{1}{q^{\prime }}, 
\]%
i.e. $(q,r)$ is nonsharp admissible. This concludes the proof of $\left(
3.10\right) $ and $(3.11)$ for nonendpoint case.

\textbf{The second step; }It remains to prove $\left( 3.10\right) $ and $%
(3.11)$ for the endpoint case, i.e. when%
\[
\left( q,\text{ }r\right) =\left( 2,\text{ }\frac{2n}{n-2}\right) ,n>2. 
\]%
It suffices to show $(3.14)$. \ By decomposing $T(F,G)$ dyadically as $%
\dsum\limits_{j}T_{j}(F,G),$ where the summation is over the integers $%
\mathbb{Z}$ and%
\begin{equation}
T_{j}\left( F,G\right) =\dint\limits_{t-2^{j-1}<s\leq t-2^{j}}\langle \left(
U\left( s\right) \right) ^{\ast }F\left( s\right) ,\left( U\left( t\right)
\right) ^{\ast }G\left( t\right) \rangle dsdt  \tag{3.19}
\end{equation}%
we see that it suffices to prove the estimate 
\begin{equation}
\dsum\limits_{j}\left\vert T_{j}(F,G)\right\vert \lesssim \left\Vert
F\right\Vert _{L_{t}^{2}L_{x}^{r^{\prime }}\left( H\right) }\left\Vert
G\right\Vert _{L_{t}^{2}L_{x}^{r^{\prime }}}\text{.}  \tag{3.20}
\end{equation}%
For this aim, before we will show the following estimate 
\begin{equation}
\left\vert T_{j}(F,G)\right\vert \lesssim 2^{-j\beta \left( a,b\right)
}\left\Vert F\right\Vert _{L_{t}^{2}L_{x}^{a^{\prime }}}\left\Vert
G\right\Vert _{L_{t}^{2}L_{x}^{b^{\prime }}}  \tag{3.21}
\end{equation}%
for all $j\in \mathbb{Z}$ and all $\left( \frac{1}{a},\frac{1}{b}\right) $\
in a neighborhood of $\left( \frac{1}{r},\frac{1}{r}\right) $. For proving $%
\left( 3.21\right) $ we will use the real interpolation of Lebesque space
and sequence spaces $l_{q}^{s}$ (see e.g $\left[ \text{15}\right] ,$ \S\ %
1.18.2 ). Indeed, by $\left[ \text{15, \S\ 1.18.4.}\right] $ we have 
\begin{equation}
\left( L_{t}^{2}L_{x}^{p_{0}},L_{t}^{2}L_{x}^{p_{1}}\right) _{\theta
,2}=L_{t}^{2}L_{x}^{p,2}  \tag{3.22}
\end{equation}%
whenever $p_{0},$ $p_{1}\in \left[ 1,\infty \right] ,$ $p_{0}\neq p_{1}$ and 
$\frac{1}{p}=\frac{1-\theta }{p_{0}}+\frac{\theta }{p_{1}}$ and $\left(
l_{\infty }^{s_{0}},l_{\infty }^{s_{1}}\right) _{\theta ,1}=l_{1}^{s}$ for $%
s_{0}$, $s_{1}\in \mathbb{R}$, $s_{0}\neq s_{1}$ and%
\[
\frac{1}{s}=\frac{1-\theta }{s_{0}}+\frac{\theta }{s_{1}}, 
\]%
where 
\[
l_{q}^{s}=\left\{ u=\left\{ u_{j}\right\} _{j-1}^{\infty },u_{j}\in \mathbb{C%
}\text{, }\left\Vert u\right\Vert _{l_{q}^{s}}=\left(
\dsum\limits_{j=1}^{\infty }2^{jsq}\left\vert u_{j}\right\vert ^{q}\right) ^{%
\frac{1}{q}}<\infty \right\} . 
\]

By $\left( 3.22\right) $ the estimate $(3.21)$ can be rewritten as 
\begin{equation}
T:L_{t}^{2}L_{x}^{a^{\prime }}\times L_{t}^{2}L_{x}^{b^{\prime }}\rightarrow
l_{\infty }^{\beta \left( a,b\right) },  \tag{3.23}
\end{equation}%
where $T=\left\{ T_{j}\right\} $ is the vector-valued bilinear operator
corresponding to the $T_{j}.$ We apply Lemma 3.2 to $\left( 3.23\right) $
with $r=1$, $p=q=2$ and arbitrary exponents $a_{0},$ $a_{1}$, $b_{0}$, $%
b_{1} $ such that 
\[
\beta \left( a_{0},b_{1}\right) =\beta \left( a_{1},b_{0}\right) \neq \beta
\left( a_{0},b_{0}\right) . 
\]

Using the real interpolation space identities we obtain%
\[
T:L_{t}^{2}L_{x}^{a^{\prime },2}\times L_{t}^{2}L_{x}^{b^{\prime
},2}\rightarrow l_{1}^{\beta \left( a,b\right) } 
\]%
for all $(a,b)$ in a neighborhood of $(r,r)$. Applying this to $a=b=r$ and
using the fact that $L^{r^{\prime }}\subset L^{r^{\prime },2}$ we obtain%
\[
T:L_{t}^{2}L_{x}^{a^{\prime },2}\times L_{t}^{2}L_{x}^{b^{\prime
},2}\rightarrow l_{1}^{0} 
\]%
which implies $\left( 3.21\right) .$

Consider the multipoint Cauchy problem for forced Schrodinger equation 
\begin{equation}
i\partial _{t}u+Lu=F,\text{ }t\in \left[ 0,T\right] ,\text{ }x\in R^{n}, 
\tag{3.24}
\end{equation}

\[
u\left( t_{0},x\right) =\varphi \left( x\right)
+\dsum\limits_{k=1}^{m}\alpha _{k}u\left( \lambda _{k},x\right) ,\text{ }%
x\in R^{n},\text{ }t_{0},\text{ }\lambda _{k}\in \left[ 0,\right. T\left.
{}\right) ,\text{ }\lambda _{k}>t_{0}\text{.} 
\]

We are now ready to state the standard Strichartz estimates:

\textbf{Lemma 3.3. }Assume the Condition 2.1 are satisfied, $\varphi \in 
\mathring{W}^{\gamma ,p}\left( R^{n}\right) $ for $\gamma \geq \frac{n}{p}$
and $p\in \left[ 1,\infty \right] $. Then problem $\left( 3.24\right) $ has
a unique generalized solution.

\textbf{Proof. }By using the Fourier transform we get from $(3.24):$%
\[
i\hat{u}_{t}\left( t,\xi \right) +L\left( \xi \right) \hat{u}\left( t,\xi
\right) =\hat{F}\left( t,\xi \right) ,\text{ } 
\]%
\begin{equation}
\hat{u}\left( 0,\xi \right) =\hat{\varphi}\left( \xi \right)
+\dsum\limits_{k=1}^{m}\alpha _{k}\hat{u}\left( \lambda _{k},\xi \right) ,%
\text{ for a.e. }\xi \in R^{n}.  \tag{3.25}
\end{equation}%
where $\hat{u}\left( t,\xi \right) $ is a Fourier transform of $u\left(
t,x\right) $ with respect to $x.$

Consider the problem%
\begin{equation}
\hat{u}_{t}\left( t,\xi \right) -iL\left( \xi \right) \hat{u}\left( t,\xi
\right) =\hat{F}\left( t,\xi \right) ,\text{ }  \tag{3.26}
\end{equation}%
\[
\hat{u}\left( 0,\xi \right) =u_{0}\left( \xi \right) ,\text{ }\xi \in R^{n},%
\text{ }t\in \left[ 0,T\right] ,\text{ } 
\]%
where $u_{0}\left( \xi \right) \in \mathbb{C}$ for $\xi \in R^{n}.$ By
Condition 2.1 and by $\left[ \text{11, \S\ 1.10, \S\ 4.1}\right] $, $%
iL\left( \xi \right) $ is a generator of a strongly continuous $C_{0}$
semigroups $U_{L}\left( t,\xi =\right) e^{itL\left( \xi \right) }$ and the
Cauchy problem $(3.26)$ has a unique solution for all $\xi \in R^{n},$
moreover, the solution can be expressed as%
\begin{equation}
\hat{u}\left( t,\xi \right) =e^{itL\left( \xi \right) }u_{0}\left( \xi
\right) +\dint\limits_{t_{0}}^{t}e^{itL\left( \xi \right) \left\vert t-\tau
\right\vert }\hat{F}\left( \tau ,\xi \right) d\tau ,\text{ }t\in \left(
0,T\right) .  \tag{3.27}
\end{equation}%
Using the formula $\left( 3.27\right) $ and the condition $\left(
3.25\right) $ we get 
\begin{equation}
u_{0}\left( \xi \right) =\hat{\varphi}\left( \xi \right)
+\dsum\limits_{k=1}^{m}\alpha _{k}U_{L}\left( \lambda _{k},\xi \right)
u_{0}\left( \xi \right) +  \tag{3.28}
\end{equation}

\[
\dsum\limits_{k=1}^{m}\alpha _{k}\dint\limits_{t_{0}}^{\lambda
_{k}}U_{L}\left( \lambda _{k}-\tau ,\xi \right) \hat{F}\left( \tau ,\xi
\right) d\tau ,\text{ }\tau \in \left( 0,T\right) . 
\]%
From $\left( 3.27\right) $ and $\left( 3.28\right) $ we obtain that the
solution of problem $\left( 3.25\right) $ can be expressed as:%
\begin{equation}
\hat{u}\left( t,\xi \right) =U_{L}\left( t,\xi \right) \hat{\varphi}\left(
\xi \right) +\dsum\limits_{k=1}^{m}\alpha _{k}U_{L}\left( \lambda _{k},\xi
\right) u_{0}\left( \xi \right) +  \tag{3.29}
\end{equation}

\[
\dsum\limits_{k=1}^{m}\alpha _{k}\dint\limits_{t_{0}}^{\lambda
_{k}}U_{L}\left( \lambda _{k}-\tau ,\xi \right) \hat{F}\left( \tau ,\xi
\right) d\tau +\dint\limits_{t_{0}}^{t}U_{L}\left( t-\tau ,\xi \right) \hat{F%
}\left( \tau ,\xi \right) d\tau ,\text{ }\tau \in \left( 0,T\right) . 
\]

Then the solution of the problem $\left( 3.24\right) $ will be expressed as
the following formula:%
\[
u\left( t,x\right) =V\left( t\right) \varphi \left( x\right)
+\dsum\limits_{k=1}^{m}\alpha _{k}V_{k}\left( t,x\right)
+\dsum\limits_{k=1}^{m}\alpha _{k}G_{k}\left( t,x\right) +G_{0}\left(
t,x\right) , 
\]%
where 
\begin{equation}
V\left( t\right) =\digamma ^{-1}\left[ U_{L}\left( t,\xi \right) \hat{\varphi%
}\left( \xi \right) \right] ,\text{ }V_{k}\left( t,x\right) =\digamma ^{-1}%
\left[ U_{L}\left( \lambda _{k},\xi \right) \hat{\varphi}\left( \xi \right) %
\right] ,  \tag{3.30}
\end{equation}%
\[
G_{k}\left( t,x\right) =\digamma ^{-1}\left[ \dint\limits_{t_{0}}^{\lambda
_{k}}U_{L}\left( \lambda _{k}-\tau ,\xi \right) \hat{F}\left( \tau ,\xi
\right) d\tau \right] ,\text{ } 
\]%
\[
G_{0}\left( t,x\right) =\digamma ^{-1}\left[ \dint\limits_{t_{0}}^{t}U_{L}%
\left( t-\tau ,\xi \right) \hat{F}\left( \tau ,\xi \right) d\tau \right] . 
\]

\textbf{Theorem 3.2}. Assume the Conditions 2.1 and 3.1 are satisfied. Let $%
\ 0\leq s\leq 1,$ $\varphi \in \mathring{W}^{s,2}\left( R^{n}\right) $, $%
F\in N^{0}\left( \left[ 0,T\right] ;\mathring{W}^{s,2}\left( R^{n}\right)
\right) $ and let $u$ : $\left[ 0,T\right] \times R^{n}\rightarrow \mathbb{C}
$ be a solution to $\left( 3.24\right) $. Then%
\begin{equation}
\left\Vert \left\vert \nabla \right\vert ^{s}u\right\Vert _{S^{0}\left( 
\left[ 0,T\right] \right) }+\left\Vert \left\vert \nabla \right\vert
^{s}u\right\Vert _{C^{0}\left( \left[ 0,T\right] ;L^{2}\left( R^{n}\right)
\right) }\lesssim   \tag{3.31}
\end{equation}%
\[
\left\Vert \left\vert \nabla \right\vert ^{s}\varphi \right\Vert
_{L^{2}\left( R^{n}\right) }+\left\Vert \left\vert \nabla \right\vert
^{s}F\right\Vert _{N^{0}\left( \left[ 0,T\right] \right) }.
\]

\textbf{Proof.} Let $2\leq q,r,\tilde{q},\tilde{r}$ $\leq \infty $ with 
\[
\frac{2}{q}+\frac{n}{r}=\frac{2}{\tilde{q}}+\frac{n}{\tilde{r}}=\frac{n}{2}. 
\]%
If $n=2,$ we also require that $q,$ $\tilde{q}>2.$ Consider first, the
nonendpoint case. By Lemma 3.3 the problem has a solution. The linear
operators in $(3.10)$ and $(3.11)$ are adjoint of one another; thus, by the
method of $TT^{\ast }$ both will follow once we prove 
\begin{equation}
\left\Vert \dint\limits_{s<t}U_{L}\left( t-s\right) F\left( s\right)
ds\right\Vert _{L_{t}^{q}L_{x}^{r}}\lesssim \left\Vert F\right\Vert
_{L_{t}^{q^{\prime }}L_{x}^{r^{\prime }}}.  \tag{3.32}
\end{equation}

Apply Theorem 3.1 with $Q=L_{x}^{2}\left( R^{n}\right) =L_{x}^{2}.$ The
energy estiamate $\left( 3.10\right) $:%
\[
\left\Vert U_{L}\left( t\right) f\right\Vert _{L_{x}^{2}}\lesssim \left\Vert
f\right\Vert _{L_{x}^{2}}
\]%
follows from Plancherel's theorem, the untruncated decay estimate%
\[
\left\Vert U_{L}\left( t-s\right) f\right\Vert _{L_{x}^{\infty }}\lesssim
\left\vert t-s\right\vert ^{-\frac{n}{2}}\left\Vert f\right\Vert
_{L_{x}^{1}},
\]%
and explicit representation of the Schr\"{o}dinger evolution operator $%
U_{L}\left( t\right) f\left( x\right) $. \ In view of $\left( 3.30\right) ,$
due to properties gropes $U_{L}\left( t\right) $ and by the dispersive
estimate $(3.4)$ we have 
\[
\left\vert \Phi \right\vert \lesssim \dint\limits_{s<t}\left\vert
U_{L}\left( t-s\right) ds\right\vert _{B\left( H\right) }\left\vert F\left(
s\right) \right\vert ds\lesssim \dint\limits_{\mathbb{R}}\left\vert
t-s\right\vert ^{-n\left( \frac{1}{2}-\frac{1}{p}\right) }\left\vert F\left(
s\right) \right\vert ds,
\]%
where 
\[
\Phi =\dint\limits_{s<t}U_{L}\left( t-s\right) F\left( s\right) ds.
\]

Moreover, from above estimate by the Hardy-Littlewood-Sobolev inequality, we
get

\begin{equation}
\left\Vert \Phi \right\Vert _{L_{t}^{q}L_{x}^{r}\left( R^{n+1}\right)
}\lesssim \left\Vert \dint\limits_{\mathbb{R}}\left\vert t-s\right\vert
^{-n\left( \frac{1}{2}-\frac{1}{p}\right) }\left\Vert F\left( s\right)
\right\Vert _{L_{x}^{r^{\prime }}\left( R^{n}\right) }ds\right\Vert
_{L_{t}^{q}\left( \mathbb{R}\right) }\lesssim  \tag{3.33}
\end{equation}%
\[
\left\Vert F\right\Vert _{L_{t}^{q_{1}}L_{x}^{r^{\prime }}}, 
\]%
where 
\[
\frac{1}{q_{1}}=\frac{1}{q}+\frac{1}{p}+\frac{1}{2}-\frac{\alpha }{n}. 
\]

The argument just presented also covers $(3.33)$ in the case $q=\tilde{q},$ $%
r=\tilde{r}$. It allows to consider the estimate in dualized form:%
\begin{equation}
\left\vert \dint \dint\limits_{s<t}\langle U_{L}\left( t-s\right) F\left(
s\right) ,G\left( t\right) \rangle ds\right\vert \lesssim \left\Vert
F\right\Vert _{L_{t}^{q^{\prime }}L_{x}^{r^{\prime }}}\left\Vert
G\right\Vert _{L_{t}^{\tilde{q}_{1}}L_{x}^{\tilde{r}^{\prime }}}  \tag{3.34}
\end{equation}%
when 
\[
\frac{1}{\tilde{q}_{1}}=\frac{1}{\tilde{q}}+\frac{1}{\tilde{p}}+\frac{1}{2}-%
\frac{\nu }{n}. 
\]%
The case $\tilde{q}=\infty ,$ $\tilde{r}=2$ follows from $(3.33)$, i.e. 
\begin{equation}
K\lesssim \left\Vert \dint\limits_{s<t}U_{L}\left( t-s\right) F\left(
s\right) ds\right\Vert _{L_{t}^{\infty }L_{x}^{2}}\left\Vert G\right\Vert
_{L_{t}^{1}L_{x}^{2}}\lesssim  \tag{3.35}
\end{equation}%
\[
\left\Vert F\right\Vert _{L_{t}^{q_{1}}L_{x}^{r^{\prime }}}\left\Vert
G\right\Vert _{L_{t}^{1}L_{x}^{2}}, 
\]%
where 
\[
K=\left\vert \dint \dint\limits_{s<t}\langle U_{L}\left( t-s\right) F\left(
s\right) ,G\left( t\right) \rangle ds\right\vert . 
\]%
From $\left( 3.3.5\right) $ we obtain the esimate $\left( 3.34\right) $ when 
$s=0.$ The general case is obtained by using the same argument.

Now, consider the endpoint case, i.e. $\left( q,r\right) =\left( 2,\frac{2n}{%
n-2}\right) $. It is suffices to show the following estimates%
\begin{equation}
\left\Vert U_{L}\left( t\right) \varphi \right\Vert
_{L_{t}^{q}L_{x}^{r}}\lesssim \left\Vert \varphi \right\Vert _{W^{s,2}\left(
R^{n}\right) },  \tag{3.36}
\end{equation}%
\begin{equation}
\left\Vert U_{L}\left( t\right) \varphi \right\Vert _{C^{0}\left(
L_{x}^{2}\right) }\lesssim \left\Vert \varphi \right\Vert _{W^{s,2}\left(
R^{n}\right) },  \tag{3.37}
\end{equation}

\begin{equation}
\left\Vert \dint\limits_{s<t}U_{L}\left( t-s\right) F\left( s\right)
ds\right\Vert _{L_{t}^{q}L_{x}^{r}}\lesssim \left\Vert F\right\Vert _{L_{t}^{%
\tilde{q}^{\prime }}L_{x}^{\tilde{r}^{\prime }}},  \tag{3.38}
\end{equation}%
\begin{equation}
\left\Vert \dint\limits_{s<t}U_{L}\left( t-s\right) F\left( s\right)
ds\right\Vert _{C^{0}L_{x}^{2}}\lesssim \left\Vert F\right\Vert
_{L_{t}^{q^{\prime }}L_{x}^{\tilde{r}^{\prime }}}.  \tag{3.39}
\end{equation}

Indeed, applying Theorem 3.1 for 
\[
Q=L^{2}\left( R^{n}\right) ,U\left( t\right) =\chi _{\left[ 0,T\right]
}U_{L}\left( t\right) 
\]%
with the energy estimate 
\[
\left\Vert U\left( t\right) f\right\Vert _{L^{2}\left( R^{n}\right)
}\lesssim \left\Vert f\right\Vert _{L^{2}\left( R^{n}\right) } 
\]%
which follows from Plancherel's theorem, the untruncated decay estimate $%
\left( 3.8\right) $ and by using of Lemma 3.1 we obtain the estimates $%
\left( 3.36\right) $ and $\left( 3.38\right) .$ Let us temporarily replace
the $C_{t}^{0}L_{x}^{2}$ norm in estimates $\left( 3.36\right) $, $\left(
3.38\right) $ by the $L_{t}^{\infty }L_{x}^{2}.$ Then, all of the above the
estimates will follow from Theorem 3.1, once we show that $U\left( t\right) $%
\ obeys the energy estimate $\left( 3.7\right) $ and the truncated decay
estimate $(3.9)$. The estimate $\left( 3.7\right) $ is obtain immediate from
Plancherel's theorem, and $\left( 3.9\right) $ follows in a similar way as
in $\left[ \text{13, p. 223-224}\right] $. To show that the operator 
\[
GF\left( t\right) =\dint\limits_{s<t}U_{L}\left( t-s\right) F\left( s\right)
ds 
\]%
is continuous in $L^{2}\left( R^{n}\right) ,$ we use the the identity 
\[
GF\left( t+\varepsilon \right) =U\left( \varepsilon \right) GF\left(
t\right) +G\left( \chi _{\left[ t,t+\varepsilon \right] }F\right) \left(
t\right) , 
\]%
the continuity of $U\left( \varepsilon \right) $ as an operator on $%
L^{2}\left( R^{n}\right) $, and the fact that 
\[
\left\Vert \chi _{\left[ t,t+\varepsilon \right] }F\right\Vert _{L_{t}^{%
\tilde{q}^{\prime }}L_{x}^{\tilde{r}^{\prime }}}\rightarrow 0\text{ as }%
\varepsilon \rightarrow 0. 
\]

From the estimates $\left( 3.36\right) -\left( 3.39\right) $ we obtain $%
\left( 3.31\right) $ for endpoint case.

\begin{center}
\textbf{4. Strichartz type estimates for solution} \textbf{to nonlinear Schr%
\"{o}dinger equation}
\end{center}

\bigskip Consider the multipoint initial-value problem $\left( 1.1\right)
-\left( 1.2\right) .$

\textbf{Condition 4.1.} Assume that the function $F:$ $\mathbb{C}\rightarrow
C$ is continuously differentiable and obeys the power type estimates

\begin{equation}
F\left( u\right) =O\left( \left\vert u\right\vert ^{1+p}\right) ,\text{ }%
F_{u}\left( u\right) =O\left( \left\vert u\right\vert ^{p}\right) ,\text{ } 
\tag{4.1}
\end{equation}

\begin{equation}
F_{u}\left( \upsilon \right) -F_{u}\left( w\right) =O\left( \left\vert
\upsilon -w\right\vert ^{\min \left\{ p,1\right\} }+\left\vert w\right\vert
^{\max \left\{ 0,p-1\right\} }\right)  \tag{4.2}
\end{equation}%
for some $p>0,$ where $F_{u}\left( u\right) $\ denotes the derivative of
operator function $F$ with respect to $u.$

From $\left( 4.1\right) $ we obtain 
\begin{equation}
\left\vert F\left( u\right) -F\left( \upsilon \right) \right\vert \lesssim
\left\vert u-\upsilon \right\vert \left( \left\vert u\right\vert
^{p}+\left\vert \upsilon \right\vert ^{p}\right) .  \tag{4.3}
\end{equation}

\bigskip \textbf{Remark 4.1. }The model example of a nonlinearity obeying
the conditions above is $F(u)=$ $\left\vert u\right\vert ^{p}u$, $p\in
\left( 1,\infty \right) $ for which the critical homogeneous Sobolev space
is $\mathring{W}_{x}^{s_{c},2}\left( R^{n}\right) $ with $s_{c}:=\frac{n}{2}-%
\frac{2}{p}.$

\textbf{Definition 4.1.} A function $F$ : $\left[ 0,T\right] \times
R^{n}\rightarrow \mathbb{C}$ is called a (strong) solution to $(1.1)-\left(
1.2\right) $ if it lies in the class 
\[
C_{t}^{0}\left( \left[ 0,T\right] ;\mathring{W}_{x}^{s,2}\left( R^{n}\right)
\right) \cap L_{t}^{p+2}L_{x}^{\frac{np\left( p+2\right) }{4}}\left( \left[
0,T\right] \times R^{n}\right) 
\]%
and obey: 
\begin{equation}
u\left( t,x\right) =V\left( t\right) \varphi \left( x\right)
+\dsum\limits_{k=1}^{m}\alpha _{k}V_{k}\left( t,x\right)
+\dsum\limits_{k=1}^{m}\alpha _{k}G_{k}\left( t,x\right) +G_{0}\left(
t,x\right) ,  \tag{4.4}
\end{equation}%
where 
\[
V\left( t\right) =\digamma ^{-1}\left[ U_{L}\left( t,\xi \right) \hat{\varphi%
}\left( \xi \right) \right] ,\text{ }V_{k}\left( t,x\right) =\digamma ^{-1}%
\left[ U_{L}\left( \lambda _{k},\xi \right) \hat{\varphi}\left( \xi \right) %
\right] , 
\]%
\begin{equation}
G_{k}\left( t,x\right) =\digamma ^{-1}\left[ \dint\limits_{t_{0}}^{\lambda
_{k}}U_{L}\left( \lambda _{k}-\tau ,\xi \right) \hat{F}\left( \tau ,\xi
\right) d\tau \right] ,\text{ }  \tag{4.5}
\end{equation}%
\[
G_{0}\left( t,x\right) =\digamma ^{-1}\left[ \dint\limits_{t_{0}}^{t}U_{L}%
\left( t-\tau ,\xi \right) \hat{F}\left( \tau ,\xi \right) d\tau \right] . 
\]
We say that $u$ is a global solution if $T=\infty $.

Let $B\left( x,\delta \right) $ denotes the boll in $R^{n}$ centered in $x$
with radius $\delta $ and $M$ denote the Hardy-Littlewood type maximal
operator that is defined as:%
\[
Mf\left( x\right) =\sup_{\delta >0}\left[ \mu \left( B\left( x,\delta
\right) \right) \right] ^{-1}\dint\limits_{B\left( x,\delta \right)
}\left\vert f\left( y\right) \right\vert dy.
\]

For proving the main result of this section we need the following:

\textbf{Proposition 4.1 }$\left[ \text{12}\right] $(Ch.2, \S\ 1, Theorem 1)%
\textbf{\ }Let $1<p<\infty $, $1<q\leq \infty .$ Then there exists a
constant $C\left( p,q\right) $ such that for all $\left\{ f\right\} _{k\geq
0}\in L^{p}\left( R^{n}\right) $ one has%
\[
\left\Vert \left\{ Mf\right\} _{k\geq 0}\right\Vert _{L^{p}\left(
R^{n}\right) }\leq C\left( p,q\right) \left\Vert \left\{ f\right\} _{k\geq
0}\right\Vert _{L^{p}\left( R^{n};l_{q}\right) }. 
\]

\textbf{Lemma 4.1 }$\left[ \text{4, Proposition 3.1 }\right] $. Assume%
\textbf{\ }$F\in C^{\left( 1\right) }\left( \mathbb{R}\right) $. Suppose $%
\alpha \in \left( 0,1\right) ,$ $1<p,$ $q,$ $r<\infty $\ \ and $%
r^{-1}=p^{-1}+q^{-1}.$ If $u\in L^{\infty }\left( \mathbb{R}\right) ,$ $%
D^{\alpha }u\in L^{q}\left( \mathbb{R}\right) $ and $F^{\prime }\left(
u\right) \in L^{p}\left( \mathbb{R}\right) $, then $D^{\alpha }\left(
F\left( u\right) \right) \in L^{r}\left( \mathbb{R}\right) $ and 
\[
\left\Vert D^{\alpha }\left( F\left( u\right) \right) \right\Vert
_{L^{r}\left( \mathbb{R}\right) }\lesssim \left\Vert F^{\prime }\left(
u\right) \right\Vert _{L^{p}\left( \mathbb{R}\right) }\left\Vert D^{\alpha
}u\right\Vert _{L^{q}\left( \mathbb{R}\right) }.
\]

\textbf{Theorem 4.1. }Assume the Cond\i tons 2.1, 3.1, 4.1 are satisfied.
Let $0\leq s\leq 1,$ $\varphi \in \mathring{W}^{s,2}\left( R^{n}\right) $
and $n\geq 1.$ Then there exists $\eta _{0}=\eta _{0}\left( n\right) >0$
such that if $0<\eta \leq \eta _{0}$ such that 
\begin{equation}
\left\Vert \left\vert \nabla \right\vert ^{s}U_{L}\left( t\right) \varphi
\right\Vert _{L_{t}^{p+2}L_{x}^{\sigma }\left( \left[ 0,T\right] \times
R^{n}\right) }\leq \eta ,  \tag{4.5}
\end{equation}%
then here exists a unique solution $u$ to $(1.1)-\left( 1.2\right) $ on $%
\left[ 0,T\right] \times R^{n}.$ Moreover, the following estimates hold

\begin{equation}
\left\Vert \left\vert \nabla \right\vert ^{s}U_{L}u\right\Vert
_{L_{t}^{p+2}L_{x}^{\sigma }\left( \left[ 0,T\right] \times R^{n}\right)
}\leq 2\eta ,  \tag{4.6}
\end{equation}%
\begin{equation}
\left\Vert \left\vert \nabla \right\vert ^{s}u\right\Vert _{S^{0}\left( 
\left[ 0,T\right] \times R^{n}\right) }+\left\Vert u\right\Vert
_{C^{0}\left( \left[ 0,T\right] ;\mathring{W}^{s,2}\left( R^{n}\right)
\right) }\lesssim \left\Vert \left\vert \nabla \right\vert ^{s}\varphi
\right\Vert _{L_{x}^{2}\left( R^{n}\right) }+\eta ^{1+p},  \tag{4.7}
\end{equation}%
\begin{equation}
\left\Vert u\right\Vert _{S^{0}\left( \left[ 0,T\right] \times
R^{n};H\right) }\lesssim \left\Vert \varphi \right\Vert _{L_{x}^{2}\left(
R^{n};H\right) },\text{ }r=r\left( p,n\right) =\frac{2n\left( p+2\right) }{%
2\left( n-2\right) +np}.  \tag{4.8}
\end{equation}

\textbf{Proof. }We apply the standard fixed point argument. More precisely,
using the Strichartz estimates $\left( 3.31\right) $, we will show that the
solution map $u\rightarrow \Phi (u)$ defined by $\left( 4.4\right) -\left(
4.5\right) $ is a contraction on the set $B_{1}\cap B_{2}$ under the metric
given by 
\[
d\left( u,\upsilon \right) =\left\Vert u-\upsilon \right\Vert
_{L_{t}^{p+2}L_{x}^{r}\left( \left[ 0,T\right] \times R^{n}\right) },
\]%
where 
\[
B_{1}=\left\{ u\in W^{\infty ,s_{c},2}=L_{t}^{\infty }W_{x}^{s_{c},2}\left( 
\left[ 0,T\right] \times R^{n}\right) :\right. 
\]

\[
\left. \left\Vert u\right\Vert _{W^{\infty ,s_{c},2}}\leq 2\left\Vert
\varphi \right\Vert _{W_{x}^{s_{c},2}\left( R^{n}\right) }+C\left( n\right)
\left( 2\eta \right) ^{1+p}\right\} , 
\]%
\[
B_{2}=\left\{ u\in W^{p+2,s_{c},r}=L_{t}^{p+2}W_{x}^{s_{c},r}\left( \left[
0,T\right] \times R^{n}\right) :\right. 
\]%
\[
\left\Vert \left\vert \nabla \right\vert ^{s_{c}}u\right\Vert
_{L_{t}^{p+2}L_{x}^{r}\left( \left[ 0,T\right] \times R^{n}\right) }\leq
2\eta \text{, }\left. \left\Vert u\right\Vert _{L_{t}^{p+2}L_{x}^{r}\left( %
\left[ 0,T\right] \times R^{n}\right) }\leq 2C\left( n\right) \text{\ }%
\left\Vert \varphi \right\Vert _{L_{x}^{2}\left( R^{n}\right) }\right\} , 
\]%
here $C(n)$ denotes the constant from the Strichartz inequality in $\left(
3.25\right) .$

Note that both $B_{1}$ and $B_{2}$ are closed in this metric. Using the
Strichartz estimate $\left( 3.31\right) $, Proposition 4.1 and Sobolev
embedding in fractional Sobolev spaces $\left[ \left[ \text{15}\right] ,%
\text{ \S }\ \text{2.3}\right] $ we get that for $u\in B_{1}\cap B_{2}$, 
\[
\left\Vert \Phi \left( u\right) \right\Vert _{L_{t}^{\infty
}W^{s_{c},2}\left( \left[ 0,T\right] \times R^{n}\right) }\leq \left\Vert
\varphi \right\Vert _{W_{x}^{s_{c},2}\left( R^{n}\right) }+C\left( n\right)
\left\Vert \langle \nabla \rangle ^{s_{c}}F\left( u\right) \right\Vert
_{_{L_{t}^{\left( p+2\right) /\left( p+1\right) }L_{x}^{r_{1}}}}\leq 
\]%
\[
\left\Vert \varphi \right\Vert _{W_{x}^{s_{c},2}\left( R^{n}\right)
}+C\left( n\right) \left\Vert \langle \nabla \rangle ^{s_{c}}u\right\Vert
_{L_{t}^{p+2}L_{x}^{\sigma }}\left\Vert u\right\Vert
_{L_{t}^{p+2}L_{x}^{np\left( p+2\right) /4}}\leq 
\]%
\[
\left\Vert \varphi \right\Vert _{W_{x}^{s_{c},2}\left( R^{n}\right)
}+C\left( n\right) \left( 2\eta +2C\left( n\right) \left\Vert \varphi
\right\Vert _{L_{x}^{2}\left( R^{n}\right) }\right) \left\Vert \left\vert
\nabla \right\vert ^{s_{c}}u\right\Vert _{L_{t}^{p+2}L_{x}^{r}}\leq 
\]%
\[
\left\Vert \varphi \right\Vert _{W_{x}^{s_{c},2}\left( R^{n}\right)
}+C\left( n\right) \left( 2\eta +2C\left( n\right) \left\Vert \varphi
\right\Vert _{L_{x}^{2}\left( R^{n}\right) }\right) \left( 2\eta \right)
^{p},
\]%
where 
\[
L_{t}^{q}L_{x}^{r}=L_{t}^{q}L_{x}^{r}\left( \left[ 0,T\right] \times
R^{n}\right) \text{, }r_{1}=r_{1}\left( p,n\right) =\frac{2n\left(
p+2\right) }{2\left( n+2\right) +np}.
\]

Similarly, 
\[
\left\Vert \Phi \left( u\right) \right\Vert _{L_{t}^{p+2}L_{x}^{r}}\leq
C\left( n\right) \left\Vert \varphi \right\Vert _{L_{x}^{2}\left(
R^{n};H\right) }+C\left( n\right) \left\Vert u\right\Vert
_{L_{t}^{p+2}L_{x}^{r}}\leq 
\]%
\[
\left\Vert \varphi \right\Vert _{W_{x}^{s_{c},2}\left( R^{n}\right)
}+2C^{2}\left( n\right) \left\Vert \varphi \right\Vert _{L_{x}^{2}\left(
R^{n}\right) }\left( 2\eta \right) ^{p}. 
\]

Arguing as above and invoking $\left( 4.5\right) ,$ we obtain 
\[
\left\Vert \left\vert \nabla \right\vert ^{s_{c}}\Phi \left( u\right)
\right\Vert _{L_{t}^{p+2}L_{x}^{r}}\leq \eta +C\left( n\right) \left\Vert
\left\vert \nabla \right\vert ^{s_{c}}\Phi \left( u\right) \right\Vert
_{L_{t}^{\left( p+2\right) /\left( p+1\right) }L_{x}^{r_{1}}}\leq 
\]%
\[
\eta +C\left( n\right) \left( 2\eta \right) ^{1+p}. 
\]

Thus, choosing $\eta _{0}=\eta _{0}\left( n\right) $ sufficiently small, we
see that for $0<$ $\eta \leq \eta _{0}$ the function $\Phi $ maps the set $%
B_{1}\cap B_{2}$ to itself. To see that it is a contraction, we repeat the
computations above and use $(4.4)$ to obtain 
\[
\left\Vert F\left( u\right) -F\left( \upsilon \right) \right\Vert
_{L_{t}^{p+2}L_{x}^{r}}\leq C\left( n\right) \left\Vert F\left( u\right)
-F\left( \upsilon \right) \right\Vert _{L_{t}^{\left( p+2\right) /\left(
p+1\right) }L_{x}^{r_{1}}}\leq 
\]%
\[
C\left( n\right) \left( 2\eta \right) ^{p}\left\Vert u-\left( \upsilon
\right) \right\Vert _{L_{t}^{p+2}L_{x}^{r}}. 
\]

Thus, choosing $\eta _{0}=\eta _{0}\left( n\right) $ small enough, we can
guarantee that is a contraction on the set $B_{1}\cap B_{2}$. By the
contraction mapping theorem, it follows that has a fixed point in $B_{1}\cap
B_{2}$. Since $\Phi $ maps into $C_{t}^{0}W_{x}^{s_{c},2}\left( \left[ 0,T%
\right] \times R^{n}\right) $ we derive that the fixed point of $\Phi $\ is
indeed a solution to $(1.1)-\left( 1.2\right) $.

In view of Definition 4.1, uniqueness follows from uniqueness in the
contraction mapping theorem.

\textbf{References}\ \ 

\begin{enumerate}
\item J. Bergh and J. Lofstrom, Interpolation spaces: An introduction,
Springer-Verlag, New York, 1976.

\item J. Bourgain, Global solutions of nonlinear Schrodinger equations.
American Mathematical, Society Colloquium Publications, 46. American
Mathematical Society, Providence, RI, 1999, MR1691575.

\item T. Cazenave and F. B. Weissler, The Cauchy problem for the critical
nonlinear Schrodinger equation in $H^{s}$. Nonlinear Anal. 14 (1990),
807-836.

\item M. Christ and M. Weinstein, Dispersion of small amplitude solutions of
the generalized Korteweg-de Vries equation. J. Funct. Anal. 100 (1991),
87-109.

\item J. Ginibre and G. Velo, Smoothing properties and retarded estimates
for some dispersive evolution equations, Comm. Math. Phys. 123 (1989),
535--573.

\item L. Escauriaza, C. E. Kenig, G. Ponce, and L. Vega, Hardy's uncertainty
principle, convexity and Schr\"{o}dinger evolutions, J. European Math. Soc.
10, 4 (2008) 883--907.

\item G. Grillakis, On nonlinear Schrodinger equations. Comm. PDE 25 (2000),
1827\{1844.

\item C. E. Kenig and F. Merle, Global well-posedness, scattering and blow
up for the energycritical, focusing, nonlinear Schrodinger equation in the
radial case. Invent. Math. 166 (2006), 645\{675.

\item M. Keel and T. Tao, Endpoint Strichartz estimates. Amer. J. Math. 120
(1998), 955-980.

\item R. Killip and M. Visan, Nonlinear Schrodinger equations at critical
regularity, Clay Mathematics Proceedings, v. 17, 2013.

\item A. Pazy, Semigroups of linear operators and applications to partial
differential equations. Springer, Berlin, 1983.

\item E. M. Stein, Singular Integrals and differentiability properties of
functions, Princeton Univ. Press, Princeton. NJ, 1970.

\item C. D. Sogge, Fourier Integrals in Classical Analysis, Cambridge
University Press, 1993.

\item R. S. Strichartz, Restriction of Fourier transform to quadratic
surfaces and deay of solutionsof wave equations. Duke Math. J. 44 (1977),
705\{714. MR0512086A.

\item H. Triebel, Interpolation theory, Function spaces, Differential
operators, North-Holland, Amsterdam, 1978.

\item T. Tao, Nonlinear dispersive equations. Local and global analysis.
CBMS Regional conference series in mathematics, 106. American Mathematical
Society, Providence, RI, 2006., MR2233925
\end{enumerate}

\end{document}